\documentclass[12pt, a4paper]{article}
\usepackage{amssymb}
\usepackage{latexsym}
\usepackage{amsbsy}

\newcommand{\bcen}{\begin{center}}     \newcommand{\ecen}{\end{center}}
\newcommand{\bay}{\begin{array}}      \newcommand{\eay}{\end{array}}
\newcommand{\beq}{\begin{eqnarray*}}      \newcommand{\eeq}{\end{eqnarray*}}
\def\az{\alpha}
\def\bz{\beta}

\def\dz{\delta}

\def\lz{\lambda}

\def\card{\mbox{card}}
\def\char{\mbox{char}}
\def\gl{\mbox{gl.dim.}}
\def\hch{\mbox{hch.dim.}}
\def\hh{\mbox{hh.dim.}}

\def\Hom{\mbox{Hom}}
\def\rank{\mbox{rank}}
\def\rad{\mbox{rad}}

\def\Ext{\mbox{Ext}}

\def\dim{\mbox{dim}}

\def\Im{\mbox{Im}}
\def\Ker{\mbox{Ker}}
\def\Coker{\mbox{Coker}}

\def\Tor{\mbox{Tor}}

\def\la{\langle}
\def\ra{\rangle}
\begin{document}

\title{\large {\bf Hochschild (co)homology dimension \footnote{Project 10201004 supported by NSFC.}}}

\author{\large Yang Han}

\date{\footnotesize Institute of Systems Science, Academy of Mathematics and Systems Science,\\ Chinese Academy
of Sciences, Beijing 100080, P.R.China.\\ E-mail:
hany@mail.iss.ac.cn}

\maketitle

\begin{abstract}

\medskip

In 1989 Happel conjectured that for a finite-dimensional algebra
$A$ over an algebraically closed field $k$, $\gl A< \infty$ if and
only if $\hch A < \infty$. Recently Buchweitz-Green-Madsen-Solberg
gave a counterexample to Happel's conjecture. They found a family
of pathological algebra $A_q$ for which $\gl A_q = \infty$ but
$\hch A_q=2$. These algebras are pathological in many aspects,
however their Hochschild homology behaviors are not pathological
any more, indeed one has $\hh A_q = \infty=\gl A_q$. This suggests
to pose a seemingly more reasonable conjecture by replacing
Hochschild cohomology dimension in Happel's conjecture with
Hochschild homology dimension: $\gl A < \infty$ if and only if
$\hh A < \infty$ if and only if $\hh A = 0$. The conjecture holds
for commutative algebras and monomial algebras. In case $A$ is a
truncated quiver algebras these conditions are equivalent to the
quiver of $A$ has no oriented cycles. Moreover, an algorithm for
computing the Hochschild homology of any monomial algebra is
provided. Thus the cyclic homology of any monomial algebra can be
read off in case the underlying field is characteristic 0.

\end{abstract}

\medskip


\bigskip

\bigskip

Let $k$ be a fixed field. All algebras considered here are of the
form $A=kQ/I$ where $Q$ is a finite quiver and $I$ is an
admissible ideal of the path algebra $kQ$. We refer to [ARS] for
the theory of quivers and their representations. It is well-known
that in case $k$ is an algebraically closed field, up to Morita
equivalence, any finite-dimensional $k$-algebra can be written as
this form (cf. [G1]).

\medskip

\bcen{\bf 1. Hochschild (co)homology}\ecen

\medskip

In this section we recall some definitions and properties on
Hochschild (co)homology.

\medskip

{\bf (1.1)} The complex $C_*(A)=(A^{\otimes *+1},b_*): \cdots
\rightarrow A^{\otimes n+1} \stackrel{b_n}{\rightarrow} A^{\otimes
n} \rightarrow \cdots \rightarrow A^{\otimes 2}
\stackrel{b_1}{\rightarrow} A \rightarrow 0$ with
$b_n(a_0,a_1,...,a_n):=
\sum^{n-1}_{i=0}(-1)^i(a_0,...,a_ia_{i+1},...,a_n)+$ \linebreak
$(-1)^n(a_na_0,a_1,...,a_{n-1}),$ where $(a_0,...,a_n)$ $:= a_0
\otimes \cdots \otimes a_n$ and $\otimes := \otimes_k$, is called
the {\it Hochschild complex} or {\it cyclic bar complex} of $A$.

The homology group $HH_n(A):=H_n(C_*(A))$ is called the $n$th {\it
Hochschild homology group} of $A$ (cf. [Lod]). Clearly $HH_n(A)$
is a $k$-vector space, we denote $\dim_kHH_n(A)$ by $hh_n(A)$.

\medskip

{\bf (1.2)} The cohomology group of the complex
$C^*(A)=(\Hom_k(A^{\otimes *},A),d^*): 0 \rightarrow A
\stackrel{d^0}{\rightarrow} \Hom_k(A,A)
\stackrel{d^1}{\rightarrow} \cdots \rightarrow \Hom_k(A^{\otimes
n},A) \stackrel{d^n}{\rightarrow} \Hom_k(A^{\otimes n+1},A)
\rightarrow \cdots$ with
$(d^n(f))(a_0,a_1,...,a_n):=a_0f(a_1,...,a_n) +
\sum^{n-1}_{i=0}(-1)^{i+1}f(a_0,...,a_ia_{i+1},...,$ \linebreak
$a_n) + (-1)^{n+1}f(a_0,...,a_{n-1})a_n,$ $HH^n(A) :=H^n(C^*(A))$,
is called the $n$th {\it Hochschild cohomology group} of $A$ (cf.
[Ho]). Clearly $HH^n(A)$ is a $k$-vector space, we denote
$\dim_kHH^n(A)$ by $hh^n(A)$.

\medskip

{\bf (1.3)} The chain complex $C'_*(A)=(A^{\otimes *+2},b'_*):
\cdots \rightarrow A^{\otimes n+2} \stackrel{b'_n}{\rightarrow}
A^{\otimes n+1} \rightarrow \cdots \rightarrow A^{\otimes 3}
\stackrel{b'_1}{\rightarrow} A^{\otimes 2} \rightarrow 0$ with
$b'_n(a_0,a_1,...,a_{n+1}):=
\sum^n_{i=0}(-1)^i(a_0,...,a_ia_{i+1},$ \linebreak $...,a_n)$ is
called the {\it bar complex} or {\it standard complex} of $A$.

The bar complex $C'_*(A)$ is a projective resolution of $A$ as an
$A^e$-module where $A^e:=A \otimes A^{op}$. Since $A \otimes_{A^e}
C'_*(A) = C_*(A)$ and $\Hom_{A^e}(C'_*(A),A)=C^*(A)$, we have
$HH_n(A)=\Tor^{A^e}_n(A,A)$ and $HH^n(A)=\Ext^n_{A^e}(A,A)$. These
two formulae are very useful to calculate Hochschild (co)homology
by constructing a convenient projective resolution of $A$ over
$A^e$.

\medskip

{\bf (1.4) Lemma.} {\it Let $A$ and $B$ be two algebras.}

(1) {\it $HH_n({\tiny \left[\begin{array}{cc} A&_AM_B\\0&B
\end{array}\right]}) = HH_n(A) \oplus HH_n(B)$ for any $A$-$B$-bimodule
$M$. In particular $HH_n(A \times B)=HH_n(A) \oplus HH_n(B)$.}

(2) {\it If $A$ and $B$ are derived equivalent then
$HH_n(A)=HH_n(B)$. In particular iterate tilting and Morita
equivalence preserves Hochschild homology.}

(3) $HH_n(A \otimes B)= \oplus _{i+j=n}HH_i(A) \otimes HH_j(B).$

\medskip

{\bf Proof.} (1) [Lod; Theorem 1.2.15]. (2) [K1; Theorem 2.2]. (3)
[Mac; Chap.VIII, Theorem 7.4]. \hfill{$\Box$}

\medskip

{\bf (1.5) Lemma.} {\it Let $A$ and $B$ be two algebras.}

(1) {\it $HH^n(A \times B)=HH^n(A) \oplus HH^n(B)$}.

(2) {\it If $A$ and $B$ are derived equivalent then
$HH^n(A)=HH^n(B)$. In particular iterate tilting and Morita
equivalence preserves Hochschild cohomology.}

(3) $HH^n(A \otimes B)= \oplus _{i+j=n}HH^i(A) \otimes HH^j(B).$

\medskip

{\bf Proof.} (1) [CE; Chap.IX, Theorem 5.3]. (2) [R; Proposition
2.5]. (3) [Mac; Chap.VIII, Theorem 7.4]. \hfill{$\Box$}

\medskip

{\bf Remark.} For Hochschild cohomology in general we have no
similar formula to Lemma (1.4) (1). An obvious counterexample is
the Kronecker algebra ${\tiny \left[\begin{array}{cc} k&k^2\\0&k
\end{array}\right]}$. Its Hochschild cohomology groups can be
obtained by [Ha; Proposition 1.6]. The relation between the
triangular matrix algebra ${\tiny \left[\begin{array}{cc}
A&_AM_B\\0&B
\end{array}\right]}$ and $A, B$ was observed in [MiP].

\medskip

\bcen{\bf 2. Hochschild (co)homology dimension}\ecen

\medskip

In this section we give the definitions and some properties of
Hochschild (co)homology dimension.

\medskip

{\bf (2.1) Definition.} ([B; Definition 3.1]) The {\it Hochschild
homology dimension} of $A$ is $\hh A := \inf \{ n \in \mathbb{N}_0
| hh_i(A)=0$ for $i > n \}$.

\medskip

{\bf Proposition.} {\it Let $A$ and $B$ be two algebras.}

(1) $\hh {\tiny \left[\begin{array}{cc} A&_AM_B\\0&B
\end{array}\right]} = \max\{\hh A, \hh B\}$ {\it for any $A$-$B$-bimodule $M$}.
{\it In particular} $\hh (A \times B)= \max\{\hh A, \hh B\}$.

(2) {\it If $A$ and $B$ are derived equivalent then} $\hh A = \hh
B$. {\it In particular iterate tilting and Morita equivalence
preserves Hochschild homology dimension.}

(3) $\hh (A \otimes B)= \hh A + \hh B$.

\medskip

{\bf Proof.} It follows from Lemma (1.4) \hfill{$\Box$}

\medskip

{\bf (2.2) Definition.} The {\it Hochschild cohomology dimension}
of $A$ is \linebreak $\hch A := \inf \{ n \in \mathbb{N}_0 |
hh^i(A)=0$ for $i > n \}$.

\medskip

{\bf Proposition.} {\it Let $A$ and $B$ be two algebras.}

(1) $\hch (A \times B)= \max\{\hch A, \hch B\}$.

(2) {\it If $A$ and $B$ are derived equivalent then} $\hch A =
\hch B$. {\it In particular iterate tilting and Morita equivalence
preserves Hochschild cohomology dimension.}

(3) $\hch (A \otimes B)= \hch A + \hch B$.

\medskip

{\bf Proof.} It follows from Lemma (1.5). \hfill{$\Box$}

\medskip

{\bf Remark.} For some special class of algebras we have some
other nice properties on Hochschild cohomology dimension:  If $A$
is hereditary then $\hch A \leq 1$, in this case $\hch A =0$ if
and only if $Q$ is a tree (cf. [Ha; Proposition 1.6 and Corollary
1.6]). If $A$ is radical square zero then $\hch A =0$ if and only
if $Q$ is a tree (cf. [Ha; Proposition 2.3]).

\medskip

\bcen{\bf 3. Relation with global dimension}\ecen

\medskip

In this section we consider the relation between the Hochschild
(co) homology dimension and the global dimension of an algebra.
This is the main part of the paper.

\medskip

{\bf (3.1) Conjecture.} ([Ha]) $\hch A < \infty$ {\it if and only
if} $\gl A < \infty$.

\medskip

Let $Q=(Q_0,Q_1)$ where $Q_0$ (resp. $Q_1$) is the set of vertices
(resp. arrows) in $Q$. For any arrow $a \in Q_1$ denote by $s(a)$
and $e(a)$ the starting point and the ending point of $a$
respectively. A path $a_1 \cdots a_n$ with $n \geq 1$ is called an
{\it oriented cycle} or a {\it cycle} of length $n$ if
$e(a_i)=s(a_{i+1})$ for $i=1,...,n-1$ and $e(a_n)=s(a_1)$.
Happel's conjecture holds for {\it truncated quiver algebras},
i.e., the algebras $kQ/k^nQ$ where $n \geq 2$ and $k^nQ$ denotes
the ideal of $kQ$ generated by all paths of length $\geq n$ (cf.
[Ci2]).

\medskip

{\bf Proposition.} ([Loc; Theorem 3]) {\it Let $A$ be a truncated
quiver algebra. Then the following are equivalent:}

(1) $\hch A < \infty$;

(2) $\gl A < \infty$;

(3) {\it The quiver of $A$ has no oriented cycles.}

\medskip

However Happel's conjecture is not true in general.

\medskip

{\bf (3.2) Pathological algebras.} Let $A:=A_q:=k \la x,y \ra /
(x^2,xy+qyx, y^2)$ with $q \in k$. The algebra $A$ is very
pathological for some special $q$. It is used to construct
counterexamples for a few problems:

(1) ([BGMS]) By constructing a minimal projective resolution of
$A$ over $A^e$ Buchweitz-Green-Madsen-Solberg show that $\hch A=2$
in the case $q$ is neither zero nor a root of unity, but $\gl A =
\infty$. This provides a negative answer to Happel's conjecture.

(2) ([Sc]) By constructing a non-projective $A_q$-module without
self-extension for some special $q$, Schulz gives a negative
answer to Tachikawa's conjecture: Over a quasi-Frobenius algebra,
a finitely generated module without self-extensions is projective.

(3) ([LS]) In case $q$ is not a root of unity, the trivial
extension of $A_q$ is a local symmetric algebra whose $AR$-quiver
contains a bounded infinite $DTr$-orbit. Thus Liu-Schulz give a
negative answer to Ringel's problem: whether the number of modules
having the same length in a connected component of the $AR$-quiver
is always finite.

(4) ([SY]) Richard showed that $A_q$ and $A_{q'}$ are socle
equivalent self-injective algebras which are not stably equivalent
in case $q' \neq q$ or $q^{-1}$.

\medskip

Though these algebras are quite pathological in many aspects,
their \linebreak Hochschild homology behaviors are not
pathological any more.

\medskip

{\bf Proposition.} {\it For any} $q \in k$, $\hh A= \infty = \gl
A$.

\medskip

{\bf Proof.} Recall the minimal projective resolution constructed
in [BGMS]:

For any $n \geq 0$ construct $\{f^n_i\}^n_{i=0} \subseteq
A^{\otimes n}:$ First $\{f^0_0\} \subseteq A^{\otimes 0} = k$ with
$f^0_0:=1$. Second $f^1_0:=x$ and $f^1_1 :=y$. Third
$f^n_{-1}:=0=:f^n_{n+1}$, and for any $n \geq 2$, $f^n_i :=
f^{n-1}_{i-1} \otimes y + q^i f^{n-1}_i \otimes x$ with $i \in \{
0,1,...,n\}$. Here $A$ is identified with $k \otimes _k A$. Note
that $f^n_i=x \otimes f^{n-1}_i + q^{n-i} y \otimes
f^{n-1}_{i-1}$.

Let $P^n:= \oplus ^n_{i=0} A \otimes_kf^n_i \otimes _k A \subseteq
A ^{\otimes (n+2)}$ for all $n \geq 0$. Set $\tilde{f}^n_i:=1
\otimes f^n_i \otimes 1$ for $n \geq 1$ and $\tilde{f}^0_0:=1
\otimes 1$. Define $\dz^n:P^n \rightarrow P^{n-1}$ by letting
$\dz^n(\tilde{f}^n_i):= x\tilde{f}^{n-1}_i+
q^{n-i}y\tilde{f}^{n-1}_{i-1}+ (-1)^n\tilde{f}^{n-1}_{i-1}y +
(-1)^nq^i\tilde{f}^{n-1}_ix$. Then $(P^*,\dz^*)$ is a minimal
projective resolution of $A$ over $A^e$.

Applying the functor $A \otimes _{A^e} -$ to the minimal
projective resolution we have $A \otimes _{A^e} (P^*,\dz^*) =
(\oplus^*_{i=0} A \otimes _k f^*_i, \tau^*)$ where $\tau^n$ is
defined by $\tau^n(\lz \otimes f^n_i)=\lz x \otimes f^{n-1}_i +
(-1)^nq^ix \lz \otimes f^{n-1}_i + q^{n-i} \lz y \otimes
f^{n-1}_{i-1} +(-1)^n y \lz \otimes f^{n-1}_{i-1}$.

Note that $A$ has a basis $\{1,y,x,yx\}$. Thus $A \otimes _k
f^n_i$ has a basis $\{ 1 \otimes f^n_i, y \otimes f^n_i, x \otimes
f^n_i, yx \otimes f^n_i \}$. Furthermore $\oplus ^n_{i=0} A
\otimes _k f^n_i$ has a basis $\{ 1 \otimes f^n_i, y \otimes
f^n_i, x \otimes f^n_i, yx \otimes f^n_i \}^n_{i=0}$. If $n \geq
1$ then $\oplus ^{n-1}_{i=0} A \otimes _k f^{n-1}_i$ has a basis
$\{ 1 \otimes f^{n-1}_i, y \otimes f^{n-1}_i, x \otimes f^{n-1}_i,
yx \otimes f^{n-1}_i \}^{n-1}_{i=0}$. Moreover $\tau^n(1 \otimes
f^n_i)= (1+(-1)^nq^i)(x \otimes f^{n-1}_i)+(q^{n-i}+(-1)^n)(y
\otimes f^{n-1}_{i-1})$, $\tau^n(y \otimes f^n_i)=
(1+(-1)^{n+1}q^{i+1})(yx \otimes f^{n-1}_i)$, $\tau^n(x \otimes
f^n_i)=(-q^{n-i+1}+(-1)^n)(yx \otimes f^{n-1}_{i-1})$ and
$\tau^n(yx \otimes f^n_i)=0$. Write $\tau^n$ as a matrix, it is
not difficult to see that $\rank(\tau^n)
\leq {\tiny \left\{ \begin{array}{ll} 2n-1,&\mbox {if $n$ odd;}\\
2n+1,&\mbox{if $n$ even.} \end{array} \right.}$ Thus
$\rank(\tau^{n+1})
\leq {\tiny \left\{ \begin{array}{ll} 2n+3,&\mbox {if $n$ odd;}\\
2n+1,&\mbox{$n$ even.} \end{array} \right.}$ Therefore
$hh_n(A)=\dim_k\Ker(\tau^n)- \dim_k\Im(\tau^{n+1}) = \dim_kP^n -
\dim_k\Im(\tau^n) - \dim_k\Im(\tau^{n+1}) = 4(n+1) - \rank(\tau^n)
- \rank(\tau^{n+1}) \geq 4(n+1) - (4n+2)=2.$ Indeed $hh_n(A) \geq
2$ for all $n$ since $hh_0(A)=\dim_k A/[A,A] =
{\tiny \left\{ \begin{array}{ll} 4,&\mbox{if $q=1$;}\\
3,&\mbox{if $q \neq 1$.} \end{array} \right.}$ \hfill{$\Box$}

\medskip

{\bf (3.3)} Proposition (3.2) suggests us to pose a seemingly more
reasonable conjecture by replacing $\hch A$ in Happel's conjecture
with $\hh A$. Indeed the following result, which is essentially
due to Keller, makes us expect more. Denote by $E$ the maximal
semisimple subalgebra of $A$, i.e., $E:=kQ_0$, and by $r$ the
Jacobson radical of $A$. Then one has $A=E \oplus r$.

\medskip

{\bf Proposition.} {\it If $\gl A < \infty$ then
$hh_i(A)={\tiny \left\{\begin{array}{ll} |Q_0|,&\mbox{\rm if }i=0;\\
0,&\mbox{\rm if }i \geq 0. \end{array} \right.}$}

\medskip

{\bf Proof.} If $\gl A < \infty$ then by [K2; Proposition 2.5] we
have \linebreak $HC_n(A, \rad A)=0$ for all $n \geq 0$. Using the
long exact sequence $$\cdots \rightarrow HH_n(A,r) \rightarrow
HC_n(A, r) \rightarrow HC_{n-2}(A,r) \rightarrow
\cdots$$$$\rightarrow HH_2(A,r) \rightarrow HC_2(A, r) \rightarrow
HC_0(A,r) \rightarrow HH_1(A,r)$$$$ \rightarrow HC_1(A, r)
\rightarrow 0 \rightarrow HH_0(A,r) \rightarrow HC_0(A, r)
\rightarrow 0$$ (cf. [Lod; E2.2.4]) we have $HH_n(A, r)=0$ for all
$n \geq 0$. Using the long exact sequence
$$\cdots \rightarrow HH_n(A,r) \rightarrow HH_n(A) \rightarrow
HH_n(E) \rightarrow $$$$ \cdots \rightarrow HH_0(A,r) \rightarrow
HH_0(A) \rightarrow HH_0(E) \rightarrow 0$$ (cf. [Lod; \S
1.1.16]), we have $HH_n(A)= HH_n(E)$ for all $n \geq 0$. Thus the
proposition follows from $hh_i(E)={\tiny \left\{\begin{array}{ll} |Q_0|&\mbox{\rm if }i=0;\\
0&\mbox{\rm if }i \geq 0. \end{array} \right.}$ (cf. [Lod;
E1.1.1]). \hfill{$\Box$}

\medskip

{\bf Corollary.} {\it If} $\gl A < \infty$ {\it then} $\hh A=0$.

\medskip

{\bf (3.4) Conjecture.} {\it The following conditions are
equivalent:}

(1) $\hh A < \infty$;

(2) $\hh A =0$;

(3) $\gl A< \infty$.

\medskip

{\bf (3.5)} Now we show that Conjecture (3.4) holds for
commutative algebras. We refer to [Mat] for commutative algebra
theory. The following result is essentially due to Avramov and
Vigue\'{e}-Poirrier.

\medskip

{\bf Theorem.} {\it Let $A$ be a commutative algebra. Then the
following conditions are equivalent:}

(1) $\hh A < \infty$;

(2) $\hh A =0$;

(3) $\gl A< \infty$

(4) $\gl A=0$.

\medskip

{\bf Proof.} By Corollary (3.1) it is enough to show that (1)
$\Rightarrow$ (4). Since $A$ is commutative, it is the direct
product of local algebras. Thus without loss of generality we may
assume that $A$ is local. If $\hh A < \infty$ then by [AV;
Theorem] $A$ is smooth. By [Mat; p.216, Lemma 1], $A$ is regular.
Furthermore $\gl A=0$ by Serre's theorem (cf. [Mat; Theorem
19.2]). \hfill{$\Box$}

\medskip

{\bf (3.6)} Now we prove that Conjecture (3.4) holds for {\it
monomial algebras} or {\it zero relation algebras} in some
literature, i.e., the algebras $A=kQ/I$ with $I$ generated by some
paths of length $\geq 2$.

A cycle $a_1 \cdots a_n$ with $\{a_1,...,a_n\} \subseteq Q_0$ is
said to be {\it basic} if the vertices $s(a_1),...,s(a_n)$ are
distinct from each other (cf. [LZ]). A cycle is said to be {\it
proper} if it is not a power ($\geq 2$) of any cycle (cf. [Ci3]).
Note that proper cycles are called basic cycles in [Sk]. Denote by
$Q^c_n$ (resp. $Q^p_n$) the set of cycles (resp. proper cycles) of
length $n$. The truncated quiver algebra whose quiver is a basic
cycle is called a {\it truncated basic cycle algebra}. The
Hochschild homology of a truncated basic cycle is known (cf. [Sk]
and [Z]). If $A$ is a monomial algebra then $A$ can be naturally
defined to be an $\mathbb{N}_0$-graded algebra by using the length
of the paths. Thus the $i$-th Hochschild homology group $HH_n(A)$
becomes an $\mathbb{N}_0$-graded $k$-vector space by
$\deg_{\mathbb{N}_0}a_0 \otimes \cdots \otimes a_n:= \sum
^n_{i=0}\deg_{\mathbb{N}_0}a_i$ for $a_i \in A, i=0,1,...,n$. The
degree $j$ part of $HH_i(A)$ is denoted by $HH_{i,j}(A)$. Let
$C_q$ be the cyclic group of order $q$. Then $C_q$ acts on $Q^c_q$
by $g \cdot (a_1 \cdots a_{q-1}a_q):=(a_qa_1 \cdots a_{q-1})$.
Under this action the orbit of the cycle $Y$ is denoted by $O_Y$.

\medskip

{\bf Proposition.} ([Sk; Theorem 2]) {\it Let $A=kQ/k^nQ$ be a
truncated quiver algebra and $q=cn+e$ for $0 \leq e \leq n-1$.
Then $HH_{p,q}(A)= $}

{\footnotesize $\left\{\begin{array}{ll} k^{a_q},&\mbox{if $1
\leq e \leq n-1$ and $2c \leq p \leq 2c+1$};\\
\oplus_{r|q}(k^{(n,r)-1} \oplus \Ker(\cdot \frac{n}{(n,r)}: k
\rightarrow
k))^{b_r},& \mbox{if $e=0$ and $0<2c-1=p$};\\
\oplus_{r|q}(k^{(n,r)-1} \oplus \Coker(\cdot \frac{n}{(n,r)}: k
\rightarrow k))^{b_r},& \mbox{if $e=0$ and $0<2c=p$;}\\
k^{|Q_0|},& \mbox{if $p=q=0$};\\ 0,& \mbox{otherwise.}
\end{array}\right.$}

{\it Here} $a_q:=\card(Q^c_q/C_q), (n,r):=\gcd(n,r)$ {\it and}
$b_r:=\card(Q^p_r/C_r)$.

\medskip

{\bf Corollary 1.} {\it Let $A=kQ/k^nQ$ where $Q$ is a basic cycle
of length $l$. Then} $hh_p(A)={\footnotesize
\left\{\begin{array}{ll}
l+[\frac{n-1}{l}],&\mbox{if $p=0$;}\\
\left[\frac{\left[\frac{p}{2}\right]n+n-1}{l}\right]-
\left[\frac{\left[\frac{p}{2}\right]n}{l}\right],&
\mbox{if $(p \geq 1) \wedge (l \nmid [\frac{p+1}{2}]n)$};\\
\left[\frac{\left[\frac{p}{2}\right]n+n-1}{l}\right]-
\left[\frac{\left[\frac{p}{2}\right]n}{l}\right]+(n,l)-1,&
\mbox{if $(p \geq 1)
\wedge (l | [\frac{p+1}{2}]n) \wedge (\char k \nmid \frac{n}{(n,l)})$};\\
\left[\frac{\left[\frac{p}{2}\right]n+n-1}{l}\right]-
\left[\frac{\left[\frac{p}{2}\right]n}{l}\right]+(n,l),& \mbox{ if
$(p \geq 1) \wedge (l | [\frac{p+1}{2}]n) \wedge (\char k |
\frac{n}{(n,l)}).$}
\end{array}\right.}$ {\it Here $[x]$ denotes the largest integer
$i$ satisfying $i \leq x$.}

\medskip

{\bf Corollary 2} {\it If $A=kQ/k^nQ$ where $Q$ is a basic cycle
of length $l$. Then} $\hh A = \infty$.

\medskip

{\bf Proof.} If $(n,l) \geq 2$ then for any $p=2ml-1$ with $m \geq
1$ one has $hh_p(A) \geq 1.$ If $(n,l)=1$ and $l \leq n-1$ then
$\left[\frac{\left[\frac{p}{2}\right]n+n-1}{l}\right] >
\left[\frac{\left[\frac{p}{2}\right]n}{l}\right]$ for any $p \geq
1$. Thus $hh_p(A) \geq 1$ for any $p \geq 1$. If $(n,l)=1$ and $l
> n-1$ then there are integers $u$ and $v$ such that $un+vl=1$.
Thus for any $p=2(ml+u-1)$ with $m \gg 0$ one has
$\left[\frac{\left[\frac{p}{2}\right]n+n-1}{l}\right]=
\left[\frac{\left[\frac{p}{2}\right]n+n-un-vl}{l}\right]=mn-v$ and
$\left[\frac{\left[\frac{p}{2}\right]n}{l}\right]< mn-v$.
Therefore $hh_p(A) \geq 1$. \hfill{$\Box$}

\medskip

Let $A=kQ/I$ be a monomial algebra. Let $C'=a_1 \cdots a_n$ be a
cycle in $Q$. Then there is a basic cycle quiver $C$ with
$|C_0|=n$ and a natural quiver morphism $C \rightarrow C'
\rightarrow Q$. We define the relations on $C$ by pulling back the
relations in $Q$, i.e., a path in $C$ is a relation if and only if
its image in $Q$ is a relation. The algebra $Z$ given by quiver
$C$ and the relations is called a {\it cycle algebra overlying
$A$}. A cycle algebra $Z$ overlying $A$ is said to be {\it
minimal} if $Z$ is not a covering of a smaller cycle algebra
overlying $A$ (cf. [IZ]). For covering theory we refer to [G2] and
[MarP]. Clearly minimal cycle algebras overlying $A$ are just
those cycle algebras overlying $A$ constructed from proper cycles
in $Q$. Denote by $HC_n(A)$ the $n$-th cyclic homology group of
$A$ (cf. [Lod]).

\medskip

{\bf Theorem 1.} {\it Let $A$ be a monomial algebra. Then $HH_n(A)
\cong \oplus _{Z}HH_n(Z)$ and $HC_n(A) \cong \oplus _{Z}HC_n(Z)$
where the index $Z$ runs through all minimal cycle algebra over
$A$.}

\medskip

{\bf Proof.} By [Ci3; Lemma 2.3] the Hochschild homology and the
cyclic homology of $A$ is just the Hochschild homology and the
cyclic homology of the $E$-{\it normalized mixed complex}
$\bar{C}_E(A)=(A \otimes_{E^e}(r^{\otimes_E n}), b, B)$ with $b$
the Hochschild boundary and $B: A \otimes _{E^e}(r^{\otimes_E n})
\rightarrow A \otimes _{E^e}(r^{\otimes_E n+1})$ given by the
formula $B(a_0,a_1,...,a_n)=
\sum^n_{i=0}(-1)^{in}(1,a_i,...,a_n,(a_0)_r,...,a_{i-1}),$ where
$a=a_E+a_r$ is the decomposition of an element $a$ of $A$ in the
direct sum $A=E \oplus r$. Note that the $k$-vector space $A
\otimes_{E^e}(r^{\otimes_E n}) = \oplus _{a_0,r_1,...,r_n} k(a_0
\otimes r_1 \otimes \cdots \otimes r_n)$ (where $a_0,r_1,...,r_n$
runs all paths sequences with $a_0r_1 \cdots r_n$ being a cycle)
$= \oplus_{\tiny \mbox{orbit $O_Y$ of proper cycle $Y$ in $Q$}}
\oplus _{a_0,r_1,...,r_n} k(a_0 \otimes r_1 \otimes \cdots \otimes
r_n)$ (where $a_0,r_1,...,r_n$ runs all paths sequences with
$a_0r_1 \cdots r_n=X^m$ for some $X \in O_Y$ and $m)=
\oplus_{\tiny \mbox{minimal cycle algebra $Z$ overlying $A$}}Z
\otimes_{F^e}(s^{\otimes_F n})$ where $F$ is the maximal
semisimple subalgebra of $Z$ and $s$ is the Jacobson radical of
$Z$. Clearly this induces an isomorphism of mixed complex from
$\bar{C}_E(A)$ to $\oplus_Z \bar{C}_F(Z)$ where $Z$ runs through
all minimal cycle algebras over $A$. By [Lod; Proposition 2.5.15]
we are done. \hfill{$\Box$}

\medskip

{\bf Algorithm.} The Hochschild homology of any monomial algebra
$A$ can be computed as following:

(1) Determine the set ${\cal O}$ of the orbits of proper cycles in
the quiver of $A$.

(2) Determine the set ${\cal M}$ of the minimal cycle algebras
overlying $A$ corresponding to the elements in ${\cal O}$.

(3) Apply the reduction of dropping projective modules in [IZ; p.
512], [IZ; Theorem 2.3] and (3.6) Corollary 1 to obtain $HH_n(Z)$
for each $Z \in {\cal M}$ and for all $n \geq 0$.

(4) Apply (3.6) Theorem 1 to obtain $HH_n(A)$ for all $n \geq 0$.

\medskip

In case $\char k=0$, the cyclic homology of $A$ can be read off:

(5) Apply [Lod; Theorem 4.1.13] to obtain $HC_n(A)$.

\medskip

{\bf Theorem 2.} {\it Let $A$ be a monomial algebra. Then the
following conditions are equivalent:}

(1) $\hh A < \infty$;

(2) $\hh A =0$;

(3) $\gl A< \infty$.

\medskip

{\bf Proof.} By (3.1) Corollary 1 it is enough to show that (1)
$\Rightarrow$ (3): If $\gl A = \infty$ then by [IZ; Proposition
3.2] there is a minimal cycle algebra $Z$ with $\gl Z = \infty$.
Via a series of dropping projective modules, $Z$ can be reduced to
a nonsemisimple truncated basic cycle algebra $Z'$ (cf. [IZ; p.
512]). By (3.6) Corollary 2 and [IZ; Theorem 2.3], one has $\hh Z=
\hh Z' = \infty$. Thus by (3.6) Theorem 1 one has $\hh A \geq \hh
Z = \infty$. It is a contradiction. \hfill{$\Box$}

\medskip

{\bf Remark 1.} The equivalence between (2) and (3) is the main
result of [IZ].

\medskip

It is well-known that the algebras whose Gabriel quiver have no
oriented cycles are of finite global dimension (cf. [ENN;
Corollary 6]) and Hochschild homology dimension 0 (cf. [Ci1]). We
show in the following that their converses hold for truncated
quiver algebras.

\medskip

{\bf Theorem 3.} {\it Let $A$ be a truncated quiver algebra (in
particular a radical square zero algebra). Then the following
conditions are equivalent:}

(1) $\hh A < \infty$;

(2) $\hh A =0$;

(3) $\gl A< \infty$;

(4) $Q$ {\it has no oriented cycles.}

\medskip

{\bf Proof.} (4) $\Rightarrow$ (3): It is well-known (cf. [ENN;
Corollary 6]). By (3.1) Corollary 1 it is enough to show that (1)
$\Rightarrow$ (4): If $Q$ has an oriented cycle then we take the
shortest cycle $Q'$ in $Q$ which must be a basic cycle. Let $I'=I
\cap kQ'$. Then $Z:=kQ'/I'$ is an truncated basic cycle algebra.
By [LZ; Theorem 2] or (3.6) Theorem 1, $HH_n(Z)$ is a direct
summand of $HH_n(A)$ for all $n \geq 0$. By (3.6) Corollary 2 we
have $\hh Z = \infty$ which implies $\hh A = \infty$. It is a
contradiction. \hfill{$\Box$}

\medskip

{\bf Remark 2.} In general the finiteness of the global dimension
does not imply that the Gabriel quiver of the algebra has no
oriented cycle. The simplest counterexample is the algebra
$A=kQ/I$ where $Q_0=\{1,2\}$, $Q_1=\{ \az : 1 \rightarrow 2, \bz :
2 \rightarrow 1 \}$ and $I= \la \bz\az \ra$. Clearly $\gl A=2$.

\end{document}